\documentclass[a4paper,12pt]{article}
\usepackage[utf8]{inputenc}
\usepackage[english]{babel}
\usepackage{amssymb, mathrsfs, amsmath, amsthm, color, enumerate, mathtools}

\usepackage[colorlinks=true, pdfstartview=FitV, linkcolor=blue, citecolor=blue, urlcolor=blue,pagebackref=false]{hyperref}

\theoremstyle{definition}
\newtheorem{problem}{Problem}
\newtheorem{theorem}{Theorem}

\newenvironment{comments}{\noindent\textit{Comments.}}{\medskip}

\newenvironment{source}{\noindent\textit{Source:}}{\medskip}
\newenvironment{relevance}{\noindent\textit{Relevance.}}{\medskip}

\newcommand{\bari}[1]{Rainger Nagel's list of problems collected in 2003 at a workshop in Bari, where the problem was posed by #1.}

\newcommand{\communicator}[1]{\noindent\textit{Communicated by} #1.\medskip}

\newenvironment{acknowledgement}{\noindent\textbf{Acknowledgement.}}{\medskip}

\begin{document}
\title{Open problems in one-parameter operator semigroups theory}

\author{Editors: Kodjo Rapha\"{e}l Madou,\footnote{McGill University, email: kodjo.madou@mail.mcgill.ca}\\
	Ivan Remizov\footnote{HSE University, email: ivremizov@yandex.ru}\, and Reihaneh Vafadar\footnote{Laval University, email: reihaneh.vafadar-seyedi-nasl.1@ulaval.ca}}
 
 
 \maketitle

	



\begin{abstract}
One-parameter strongly continuous semigroups of linear bounded operators on Banach spaces (also known as $C_0$-semigroups) are a fundamental operator-theoretic tool used in the study of linear and non-linear evolution PDEs arising in physics, probability, control theory and other areas of science and technology, including quantum theory, transportation problems and finance. Since 2021 the annual online conference on one-parameter semigroups of operators (OPSO) offers an opportunity for researchers worldwide to discuss the current state of the art and exchange knowledge. On every OPSO conference open problems are proposed, collected and discussed, some of them are found to be solved and hence get crossed out from the list. In this paper we provide a collection of open problems of semigroup theory that were proposed by participants of OPSO 2021--2024 conferences and have not been solved, yet. For each problem some comments on the relevance and and history of the problem are provided.  
\end{abstract}


\medskip


\newpage
\tableofcontents

\section{Generation problems}

\begin{problem}[Inverse generator problem]
	Let $H$ be a Hilbert space over $\mathbb{C}$ 
	and let $A: D(A) \subseteq H \to H$ be the infinitisimal generator of a bounded semigroup on $H$. Assume further that $A^{-1}$ exists as a densely defined, closed operator.
	
	Is $A^{-1}$ the infinitesimal generator of a bounded semigroup?
\end{problem}

\begin{comments}
This problem was originally posed (for Banach spaces) by R.\ de Laubenfels \cite{deLa88}. However, it is not hard to show that for general Banach spaces the answer to the problem is negative, see e.g.\ \cite{Zwar07}. The counter-example can be chosen such that the strongly continuous semigroup generated by $A$ is a contraction semigroup, whereas in Hilbert spaces the answer to the problem is positive for a generator of a contraction semigroup (almost trivial to show). 

There is a strong relation between the inverse generator problem and the question whether the Cayley transform of $A$ is power bounded, i.e.\ if $\sup_n \|A_d^n\| <\infty$, where $A_d= (I+A)(I- A)^{-1}$, see \cite{GuZw06}. The latter question is related to numerical analysis, since this Cayley transform pops up when applying the Crank-Nicolson scheme to the differential equation $\dot{x}(t) = A x(t)$. 

When the answer to the inverse generator problem is positive, then the (strong) stability of the semigroup generated by $A$ is equivalent to the (strong) stability of the semigroup generated by $A^{-1}$. Furthermore, it is equivalent to the strong stability of $A_d$, i.e., $\lim_{n\rightarrow \infty} A_d^n x = 0$ for all $x \in H$, \cite{GuZw06}. 

For finite-dimensional Hilbert spaces $H$, it is clear that the problem has a positive answer. For these spaces the question is; if there exists a constant $c$ independent of the dimension of $H$ such that $\sup_{t\geq 0} \|e^{A^{-1} t}\| \leq c \sup_{t\geq 0} \|e^{A t}\|$.

In 2017, a nice survey on the problem appeared \cite{Gomi17}. In that paper, the interested reader can find more results and references on the inverse generator problem. 
\end{comments}

\communicator{Hans Zwart}


\begin{problem}[An analytic semigroup generation problem]
This problem, arising in the theory of Gaussian open quantum systems, consists in proving that a dissipative operator $G$ 
which is quadratic in creation and annihilation operators (or, after a unitary transformation, a differential operator 
quadratic in partial derivatives $\partial_j$ and multiplication by coordinates $x_k$) generate an analytic semigroup.
The solution has applications in the proof of strong positivity, irreducibility and regularity properties of Gaussian 
open quantum systems. 
 
Let $a_j, a^*_k$ be the annihilation and creation operators on $\ell^2(\mathbb{N}^d;\mathbb{C})$ that are defined
by closure from their action on the canonical orthonormal basis $(e(n_1,\ldots,n_d))_{n\in\mathbb{N}^d}$, which is given by
\begin{eqnarray*}
a_j\,e(n_1,\ldots,n_d)& =& \sqrt{n_j}\ e(n_1,\ldots,n_{j-1},n_j-1,\ldots,n_d),\\
a_k^*\,e(n_1,\ldots,n_d)&=&\sqrt{n_k+1}\ e(n_1,\ldots,n_{k-1},n_k+1,\ldots,n_d).
\end{eqnarray*}
Let $H,L_\ell$ be the closures of operators defined on the canonical basis by
\begin{eqnarray*}
H &=& \sum_{j,k=1}^d \left( \Omega_{jk} a_j^* a_k + \frac{\kappa_{jk}}{2} a_j^* a_k^*
+ \frac{\overline{\kappa}_{jk}}{2} a_j a_k \right) + \sum_{j=1}^d \left( \frac{\zeta_j}{2}a_j^*
+ \frac{\bar{\zeta_j}}{2} a_j \right), \\
L_\ell &=& \sum_{k=1}^d \left( \overline{v}_{\ell k} a_k + u_{\ell k}a_k^*\right) 
\qquad \ell =1,\dots, 2d,
\end{eqnarray*}
where $\Omega:=(\Omega_{jk})_{1\leq j,k\leq d} = \Omega^*$ and $\kappa:= (\kappa_{jk})_{1\leq j,k\leq d}
= \kappa^{\hbox{\scriptsize T}} \in M_d(\mathbb{C})$ are $d\times d$ complex matrices with $\Omega$ Hermitian 
and $\kappa$ symmetric, $V=(v_{\ell k})_{1\leq \ell\leq 2d, 1\leq  k\leq d}$, $U=(u_{\ell k})_{1\leq \ell\leq 2d, 1\leq  k\leq d}$
are $2d\times d$ complex matrices $\zeta=(\zeta_j)_{1\leq j\leq d}\in\mathbb{C}^d$.

It is not difficult to show (see e.g.\ \cite[Proposition 4.9]{Fa}) that the closure of the operator $G$ defined 
on the canonical basis by
\[
G=-\mathrm{i}H -\frac{1}{2}\sum_{\ell=1}^{2d} L_\ell^*L_\ell
\]
generates a $C_0$ contraction semigroup $P=(P_t)_{t\geq 0}$ on $\ell^2(\mathbb{N}^d;\mathbb{C})$.

Suppose that the non-degeneracy condition (block-matrix form)
\[
\mathbb{K}=
\left[\begin{array}{c} V^{\hbox{\scriptsize T}} \\ U^*\end{array} \right] \left[\begin{array}{cc} \overline{V} & U\end{array} \right] =
\left[\begin{array}{cc} V^{\hbox{\scriptsize T}}\overline{V} & V^{\hbox{\scriptsize T}}U \\ U^*\overline{V} & U^*U \end{array} \right] >0
\]
holds, then sufficient conditions for $P$ to be an analytic semigroup are also available. 

The problem it to find a \emph{classification} of the set of parameters $\Omega,\kappa,U,V, \zeta$ 
for which $P$ is analytic. 

\smallskip
Alternatively, by the unitary correspondence of the above basis with multidimensional Hermite polynomials
(multiplied by $\exp(-|x|^2/2)$ normalized), one can formulate the problem with differential operators
\[
a_j = \frac{1}{\sqrt{2}}\left(x_j + \frac{\partial}{\partial x_j} \right), \qquad
a_k^* = \frac{1}{\sqrt{2}}\left(x_k - \frac{\partial}{\partial x_k} \right).
\]
In this case strict positivity of $\mathbb{K}$ implies 
\[
\left[\begin{array}{cc} 1 & -1 \\ 1 & 1\end{array} \right] 
\mathbb{K}\left[\begin{array}{cc} 1 & 1 \\ -1 & 1\end{array} \right] >0
\]
from which one finds the strong ellipticity condition for the self-adjoint part $G_0=-(1/2)\sum_{\ell=1}^{2d}L_\ell^*L_\ell$ 
of $G$
\[
\operatorname{Re}\sum_{j,k=1}^d
\left(U^*U+V^{\hbox{\scriptsize T}}\overline{V}-V^{\hbox{\scriptsize T}}U-U^*\overline{V}\right)_{jk}\overline{z}_jz_k 
> \epsilon \Vert z \Vert^2
\]
for $z=(z_j)_{1\leq j\leq d}\in\mathbb{C}^d$.

\smallskip
Thinking of spectra of $G_0$ and $H$ one wonders if the semigroup $P$ generated by $P$ is analytic when 
$\mathbb{K}>0$ and $H$ is bounded from below or from above.
\end{problem}

\communicator{Franco Fagnola}


\begin{problem}[Lumer--Phillips for transition semigroups]
	Let $\Omega$ be a Polish space. We call a semigroup $T= (T(t))_{t\geq 0}$ of  contractions on $C_b(\Omega)$ 
	a \emph{transition semigroup} if 
	\begin{enumerate}[(i)]
	\item  for every $f\in C_b(E)$ we have $T(t)f \to T(s)f$ locally uniformly  whenever $t\to s$ and
	\item  for every $t\geq 0$ and every uniformly bounded sequence $(f_n)_{n\in \mathbb{N}}\subset C_b(\Omega)$ 
	that converges locally uniformly to $f$, we have that $T(t)f_n$ converges locally uniformly to $f$.
	\end{enumerate}
	Is there a characterization of the generators of such a semigroup akin to the Lumer--Phillips theorem?	
\end{problem}

\begin{comments}
The name `transition semigroup' is inspired by applications in probability theory, where semigroups with the above properties frequently appear as transition semigroups of Markov processes (typically, these semigroups are additionally positive). We should point out that semigroups of this kind (at least at first glance) do not fit into the theory of semigroups on locally convex spaces (see e.g.\ \cite{yosida}) as this theory requires equicontinuity of the operators involved. But even the simplest example of the heat semigroup on $C_b(\mathbb{R}^d)$ shows that one cannot expect equicontinuity with respect to the topology $\tau_{\mathrm{co}}$ 
of uniform convergence on compact subsets of $\mathbb{R}^d$.

Consequently, in the literature several approaches were developed where this equicontinuity condition was weakened such as the `theory of of weakly continuous semigroups' by Cerrai \cite{cerrai} (where instead of $C_b(\Omega)$ one works on the space $BUC(\Omega)$) or the `theory of bi-continuous semigroups' by K\"uhnemund \cite{kuehnemund}; both approaches allow for a Hille--Yosida type generation result. On the other hand, \cite[Theorem 4.4]{kunze} shows that conditions (i) and (ii) the above definition already entail equicontinuity:  Not with respect to $\tau_{\mathrm{co}}$ but with respect to the so-called \emph{strict topology} $\beta_0$ (which agrees with $\tau_{\mathrm{co}}$ on $\|\cdot\|_\infty$-bounded subsets of $C_b(\Omega)$). This allows us to use the results from \cite{yosida} after all to characterize the generators of transition semigroups.\smallskip

Thus, characterizations of generators of transition semigroups are available in the literature. However, to the best of my knowledge, none of these Hille--Yosida type theorems was ever used to establish that a certain operator generates a transition semigroup (even though many examples of such semigroups and also their generators are known). This is not as surprising as it might seem, for even in the setting of strongly continuous semigroups the Hille--Yosida theorem is difficult to apply. This is due to the fact that this result requires us (in the case of bounded semigroups) to prove uniform boundedness of the family $\{\lambda^nR(\lambda, A)^n : \lambda >0, n\in \mathbb{N}\}$, which is difficult in concrete examples. In the case of non-strongly continuous semigroups, one would have to establish equicontinuity of this family of operators -- an even harder task.

The `weapon of choice' to prove that a given operator generates a strongly continuous semigroup is rather the Lumer--Phillips theorem (see \cite[Theorem II.3.15]{en}) which, however, only characterizes the generators of \emph{contraction} semigroups. The main advantage of the Lumer-Phillips theorem is that one does not have to consider the powers of the resolvent. Indeed, given dissipativity, we only need to check the so-called \emph{range condition}, i.e.\ we need to prove that $\lambda - A$ has a dense range for some $\lambda >0$.\smallskip

It would be very interesting to have a Lumer--Phillips type result for transition semigroups. Part of the problem is to find out how such a result should look like. Here, usability is (in my opinion) more important than generality. If we can obtain a sufficient condition for generation (which might make use of additional properties that one has at hand in many possible applications such as positivity of the resolvent or the strong Feller property for the resolvent) this would be already be very nice.
\end{comments}

\communicator{Markus Kunze}


\begin{problem}[Sum of squares]
    Let $A$ and $B$ be generators of $C_0$ semigroups. Under which condition does 
    \begin{center}
        $C:=\overline{A^{2}+B^{2}}$
    \end{center} 
    generate an analytic semigroup?
\end{problem}

\begin{source}
	\bari{Alessandra Lunardi}
\end{source}


\begin{problem}[Existence of a generator for analytic semigroups]
    Let $(T(t))_{t\geq0}$ be a not necessary strongly continuous semigroup on a Banach space $X$ and consider the following condition.	
	\begin{enumerate}[(i)]
	    \item $(0, \infty) \ni t  \mapsto T(t) \in L(X) $ is analytic.
	    \item $t \mapsto T(t) $ is analytic on a sector containing $\mathbb{R}_{+}$
	    \item $\|T(t)\| \leq Me^{t\omega}, \frac{d}{dt}T(t) \in L(X)$ and $\|\frac{d}{dt}T(t)\| \leq \frac{M}{t}e^{t\omega}$ for some constants $\omega, M$.
	\end{enumerate}
	Under which assumption added to (i), (ii) or (iii) does there exist a sectorial operator A generating $(T(t))_{t\geq0}$?
\end{problem}

\begin{source}
    \bari{Alessandra Lunardi}
\end{source}


\begin{problem}[Generator of product semigroups]
    Let $A$ and $B$ be the generators of two communicating $C_0$-semigroups on a Banach space and let $G$ be the generator of the corresponding product semigroup. Find (the most general) conditions implying $$D(G) = D(A) \cap D(B).$$
\end{problem}

\begin{comments}
    This yields abstract "maximal regularity" results.
\end{comments}

\begin{source}
    \bari{Rainer Nagel}
\end{source}


\begin{problem}[$A+B$ Problem]
	Let $A$ and $B$ be generators of $C_0$-semigroups on a Banach space $X$.

	\begin{enumerate}[(i)]
   		\item Define the sum $C=A+B$ such that $C$ becomes a (maximal) closed operator on $X$ and $Cx=\tilde{A}x+\tilde{B}x$ for all $x\in X$ and some extrapolated operators $A$ and $B$.
    	
    		\item Find assumptions on $A$ and $B$ such that $C$ remains a generator on $X$, thereby unifying known perturbation results.
	\end{enumerate}
\end{problem}

\begin{comments} 
	A test case is the following. Let $\mathcal{X} = C_0 (\mathbb{R},X)$ or $\mathcal{X} = L^{p}(\mathbb{R},X)$, and take $Af = f'$ and $Bf(s)= C(s)f(s)$ for appropriate $f \in \mathcal{X}$ and closed operator $C(s)$ on $X$.
\end{comments}

\begin{source}
    R. Nagel's list of problems was collected in 2003 at the workshop in Bari.
\end{source}


\begin{problem}[Non-autonomous abstract Cauchy problems]
	For unbounded linear operators $A(t)$ on a Banach space $X$ and for a starting time $t_0$, characterize the well-posedness of the non-autonomous Cauchy problem
	\begin{align*}
	    \dot x(t)   & = A(t) x(t) \quad \text{for } t \geq t_0 \\
	         x(t_0) & = x_0
	\end{align*}
	by Hille-Yosida type condition for an operator $G$ generating an evolution semigroup on
	$\mathcal{X} = C_0 (\mathbb{R},X)$ or $\mathcal{X} = L^{p}(\mathbb{R},X)$.
\end{problem}

\begin{source}
    R. Nagel's list of problems was collected in 2003 at the workshop in Bari.
\end{source}

\section{Regularity}


\begin{problem}[Characterization of immediate norm continuity by the resolvent]
    Prove or disprove that a semigroup with generator $A$ is immediately norm continuous if and only if
	\[
		\|R(is,A)\| \to 0 \text{ as } |s|\to\infty.
	\]
\end{problem}

\noindent \textit{Comment} (Charles Batty). This problem  was answered negatively, by Ralph Chill and Yuri Tomilov in \cite{chTo2009}, and independently by Tam\'as M\'atrai in \cite{Ma2008}.

\medskip

\begin{source}
    \bari{Charles J. K. Batty and Klaus-Jochen Engel}
\end{source}


\begin{problem}[Regulariziation property of non-$C_0$-semigroups]
    Consider "non-$C_0$-semigroups", e.g., bi-continuous semigroups, and describe appropriate regularization properties. 
    
\end{problem}

\begin{comments}
    Compare the Ornstein-Uhlebeck semigroup in $C_{b}(\mathbb{R}^{n})$, see \cite{lunardi2009}.
\end{comments}

\begin{source}
    \bari{Alessandra Lunardi}
\end{source}


\begin{problem}[Asymptotic analyticity]
    Which $C_0$-semigroups are "asymptotically analytic"?
\end{problem}

\textit{Comment.} First, we explain the notion of "asymptotically analytic semigroups". Let $B$ be a positive selfadjoint operator in Hilbert space $X$, and let $a>0$. Let $u$ be a solution of the graph equation
\begin{align*}
    u''+2au'+Bu = 0.
\end{align*}
This problem is governed by a $C_0$-semigroup on energy space based on $X$. Eckstein-Goldstein-Leggas in  \cite{EGL1999} proved that
\begin{align*}
    u(t) = v(t)+w(t),
\end{align*}
where $u$ satisfies the heat equation
\begin{align*}
    2av'+Bv = 0
\end{align*}
and $\| w(t)\| = o(\|v(t)\|)$ as $t \rightarrow \infty$. This leads to the notion of "asymptotically analytic".

Let  $(T(t))_{t\geq0}$ be a $C_0$-semigroup on a Banach space $X$, let $S:=(s(t))_{t\geq0}$ be an analytic $C_0$-semigroup on a Banach space $Y$ and let be $P$ a (somehow natural) bounded linear operator from $X$ to $Y$.
We call $S$  "asymptotically analytic" if there exist $S$, $P$ as above such that from all $f \in X$ there is $g \in Y$ so that
\begin{align*}
    T(t)f = S(t)g+w(t), \quad t\geq 0,
\end{align*}
where $\|w(t)\| = o(\|S(t)g\|)$  as $t \rightarrow \infty$.

Thus asymptotically analytic semigroups have the asymptotics of analytic semigroups, except for errors that are relatively small asymptotically.
As a first step, it would be of interest to have some results in the case of $Y = X$ and $P = I$. More generally, which $C_0$-groups (such as those governing second-order equations) and which non-analytic semigroups governing FDEs are asymptotically analytic?


\medskip

\begin{source}
	\bari{Jerome A.\ Goldstein}
\end{source}


\section{Long-term behaviour of semigroups}

\begin{problem}[Tauberian Theorem for Semigroups of Kernel Operators]
	Let $E\coloneqq L^p(\Omega,\mu)$ for $1\leq p \leq \infty$ and a $\sigma$-finite measure space $(\Omega,\mu)$
	and let $(T_t)_{t\in [0,\infty)} \subseteq \mathcal{L}(E)$ be a strongly continuous semigroup on $E$ such that $T_t$ is a positive \emph{kernel operators} for each $t>0$,
	meaning that there exists a measurable function $k_t \colon \Omega  \times \Omega \to \mathbb{R}_+$ such that

	\[ (T_t f)(y) = \int_\Omega k_t(y,x) f(x) \mathrm{d}\mu(x) \]
	for almost every $y\in \Omega$.

	Assume that $(T_t)_{t\in [0,\infty)}$ is \emph{mean ergodic}, meaning that the limit
\[	\lim_{T\to \infty} \frac{1}{T} \int_0^T T_t f \mathrm{d} t \]
exists in $E$ for all $f\in E$. 
	Does it follow that $(T_t)$ is \emph{strongly convergent}, i.e.\ $\lim_{t\to\infty} T_t f$ exists for all $f\in E$?
\end{problem}

\begin{comments}
	It is well known that positive semigroups of kernel operators are strongly convergent provided that the semigroups possess a fixed point $f\in E$ satisfying
	$f(x) >0$ for almost every $x\in \Omega$. This has been proven in a very general setting in \cite[Theorem 3.5]{gg2019}. The existence of such an "quasi-interior" fixed
	point is crucial and cannot be omitted: For instance, the Gaussian semigroup on $L^1(\mathbb{R})$ is not strongly convergent and fulfils all assumptions of \cite[Theorem 3.5]{gg2019} except 
	that it does not possess a quasi-interior fixed point.
	On the other hand, the Gaussian semigroup is not mean ergodic on $L^1(\mathbb{R})$. 
	In fact, using \cite[Theorem 3.5]{gg2019} it is not difficult to show the following Tauberian theorem:

	\begin{theorem}
		Let $(T_t)_{t \in [0,\infty)}$ be a positive, bounded and mean ergodic $C_0$-semigroup on $L^1(\Omega)$ for any measure space $\Omega$.
		If $T_{t_0}$ is kernel operator for some $t_0 > 0$, then $(T_t)_{t\in [0,\infty)}$ is strongly convergent.	
	\end{theorem}

	The theorem above and similar results for semigroups on spaces of measures can be found in \cite[Theorems 2.1, 4.1, 5.4]{ggk2020}.
	This gives rise to the conjecture that the existence of a quasi-interior fixed point 
	in \cite[Theorem 3.5]{gg2019} can always, i.e.\ for every $1\leq p\leq \infty$ or  more generally for every Banach lattice $E$,
	be omitted in case of a mean ergodic $C_0$-semigroup $(T_t)_{t \in (0,\infty)}$.
	In support of this conjecture, one may note that the Gaussian semigroup on $L^p(\mathbb{R})$ for $p \in (1,\infty)$ is mean ergodic and converges strongly to $0$. 
\end{comments}

\communicator{Moritz Gerlach}


\begin{problem}[Growth bounds]
    Prove or disprove that the non-analytic growth bound $\zeta(T)$ of a $C_0$-semigroup $T$ coincides with the critical growth bound $\omega_{crit}(T)$.
\end{problem}

\begin{comments}
    It is known that:
    \begin{align*}
        \zeta(T) \leq \omega_{ess}(T) 
        \quad \text{and} \quad
        \zeta(T) = \omega_{crit}(T)
    \end{align*}
    in each of the following cases:
    \begin{itemize}
        \item 
        $T$ is a $C_0$-semigroup on a Hilbert space,
        
        \item 
        $T$ is a $C_0$-semigroup on a Banach space, and $T$ has an $L^{p}$-resolvent for some $p \in (1,\infty)$,
        
        \item 
        $T$ is eventually differentiable.
    \end{itemize}
\end{comments}



\begin{relevance}
	This problem is interesting because $\zeta(T)$ and $\omega_{crit}(T)$ are both candidates to be variants of the exponential growth bound $\omega_0(T)$, modulo analytic functions and spectral bounds of the generator modulo horizontal strips. 
	Some variants of standard results have been obtained using $\zeta(T)$ instead of $\omega_(T)$ have been established, see \cite[Section 5]{Batty2007}. 
\end{relevance}

\begin{source}
    \bari{Charles J. K. Batty}
\end{source}


\begin{problem}[Precise growth behaviour in terms of the growth bound]
    Let $(T(t))_{t\geq0}$ be a $C_0$-semigroup with growth bound
	$$
		\omega_0 := \inf \{\omega \in \mathbb{R}:\|T(t)\| \leq M^{\omega} \cdot e^{t\omega}  \textrm{ for } t\geq0\}   
	$$
	Find condition such that $\omega_0$ is a minimum, i.e.,
 	$$\|T(t)\| \leq M_0 \cdot e^{t\omega_{0}}\textrm{ for } t \geq 0$$.
 \end{problem}

\begin{comments}
    This corresponds to a characterization of boundedness for semigroups.
\end{comments}

\begin{source}
    \bari{Rainer Nagel}
\end{source}


\begin{problem}[Fine scale of the growth behaviour]
    Let $(T(t))_{t\ge0}$ be a bounded $C_0$-semigroup on a (complex) Banach space $X$, and let $A$ denote its generator. Suppose that $\sigma(A)\cap i\mathbb{R}=\emptyset$ and that there exists $\alpha>0$ such that $\|(is-A)^{-1}\|=O(|s|^\alpha)$ as $|s|\to\infty.$ Find $\beta\in[0,1]$ depending on the geometric properties of the space $X$ (for instance its Fourier type, its type or its cotype) such that
    $$\|T(t)A^{-1}\|=O\bigg(\frac{\log(t)^{\beta/\alpha}}{t^{1/\alpha}}\bigg),\qquad t\to\infty.$$
\end{problem}

\begin{comments}
    It was shown by Batty and Duyckaerts \cite{Batty2008} that one may always take $\beta=1$; see also Chill and Seifert \cite{Chill2016}. Batty and Duyckaerts \cite{Batty2008} moreover showed that negative values of $\beta$ are in general not permissible. It is reasonable, therefore, to restrict attention to values of  $\beta$ lying in $[0,1]$. Borichev and Tomilov \cite{Borichev2010} showed that one may take $\beta=0$, yielding the best possible upper bound, if $X$ is a Hilbert space. This result may be viewed as a special case of more recent results appearing in papers by Batty, Chill and Tomilov \cite{Batty2016} and Rozendaal, Seifert and Stahn \cite{Rozendaal2019}. Borichev and Tomilov \cite{Borichev2010}  also showed, by considering the left-shift semigroup on a certain subspace of $\mathrm{BUC}(\mathbb{R}_+)$ with an appropriate norm, that one may have 
    $$\limsup_{t\to\infty}\frac{t^{1/\alpha}}{\log(t)^{1/\alpha}}\|T(t)A^{-1}\|>0.$$
    Hence one cannot in general hope to do better than $\beta=1$ unless one imposes additional assumptions on $X$; see also Debruyne and Seifert \cite{Debruyne2019a}. 
\end{comments}

\begin{relevance}
	This problem is important from a theoretical point of view, as its solution would elegantly complement our current understanding of polynomial stability of $C_0$-semigroups. Furthermore, there are likely to be interesting applications to concrete evolution equations on $L^p$-spaces and other (non-Hilbertian) Banach spaces with non-trivial geometric properties.
\end{relevance}

\communicator{David Seifert}


\section{Positivity of semigroups}

\begin{problem}[Infinite speed of propagation]
	Let $E\coloneqq L^p(\Omega,\mu)$ for $1\leq p < \infty$ and a $\sigma$-finite measure space $(\Omega,\mu)$
	and let $\coloneqq (T_t)_{t\in [0,\infty)} \subseteq \mathcal{L}(E)$ be a strongly continuous semigroup of positive kernel operators on $E$,
	meaning that for each $t>0$ there exists a measurable function $k_t \colon \Omega  \times \Omega \to \mathbb{R}_+$ such that

	\[ (T_t f)(y) = \int_\Omega k_t(y,x) f(x) \mathrm{d}\mu(x) \]
	for almost every $y\in \Omega$.

	Assume that $(T_t)$ is \emph{irreducible}, i.e.\ for all non-zero and positive $f,g \in E_+$ there exists $t\geq 0$ such that $T_t f \wedge g \neq 0$ (where $\wedge$ denotes the infimum
	in the Banach lattice $E$).  Does it follow that $(T_t)$ is \emph{expanding}, i.e.\ $T_tf>0$ almost everywhere for all $t>0$ and all non-zero positive $f\in E_+$?
\end{problem}

\begin{comments}
The notion of irreducibility is of great importance, for instance in the spectral theory of positive semigroups and obviously, every expanding semigroup is irreducible. On the other hand, the rotation semigroup on the unit circle is an easy example of an irreducible semigroup that fails to be expanding.

However, under certain circumstances the two properties, irreducible and expanding, are equivalent.
For instance, every holomorphic positive semigroup is known to be expanding if it is irreducible \cite[Theorem C-III 3.2]{nagel1986}.
A surprisingly little-known fact is that the same holds for all positive semigroups on atomic Banach lattices like $\ell^p$ for $1\leq p <\infty$. 
This is due to a fact which is referred to as ``L\'evy's Theorem'' in literature:
for every stochastic transition matrix $(p_{i,j})$ one either has $p_{i,j}(t) = 0$ or $p_{i,j}(t) > 0$ for all $t>0$.
The proof given by K.\ L.\ Chung in the appendix of \cite{chung1963} for this statement can easily be transferred to the setting of semigroups on atomic spaces.
Since operators on atomic spaces often serve as a prototype for general kernel operators, it is  natural to ask whether the same implication is true for semigroups of kernel operators.

The notation \emph{expanding} is not used in a uniform manner in the literature; the property is, for instance, also referred to as ``strongly positive'' or ``positivity improving''.
We would like to advertise the more systematic naming convention of \cite[Section 9.1]{abramovich2002}.
\end{comments}

\noindent \textit{Comment} (Jochen Glück). The answer to the problem is negative, here is a counterexample: 
for each $t > 0$ let $k_t: \mathbb{R} \to [0,\infty)$ by given by
\begin{align*}
	k_t(x) = 
	\begin{cases}
		\frac{1}{\Gamma(t)} x^{t-1} e^{-t} \quad & \text{if } x > 0, \\ 
		0 \quad & \text{if } x \le 0,
	\end{cases}
\end{align*}
where $\Gamma$ denotes the Gamma function. 
So $k_t$ is the probability density function of a Gamma distribution with expected value $t$ and variance $t$. 
Then $(k_t)_{t \in (0,\infty)}$ is a semigroup with respect to convolution 
and convolution with $k_t$ defines a $C_0$-semigroup $(T(t))_{t \ge 0}$ on $L^1(\mathbb{R})$. 
Now let $(S(t))_{t \ge 0}$ denote the left shift semigroup on $L^1(\mathbb{R})$. 
Then $(T(t)S(t))_{t \ge 0}$ is a positive and irreducible $C_0$-semigroup on $L^1(\mathbb{R})$ that consists of kernel operators, 
but is not expanding.

\medskip 

\communicator{Moritz Gerlach} 


\begin{problem}[Positive commutator problem]
	Let $E$ be a Banach lattice	and let $C\colon E \to E$ be a positive quasinilpotent compact operator.
    Do there exist positive operators, $A,B\colon E \to E$ such that $C=AB-BA$ with one of $A$ and $B$ compact?
\end{problem}

\begin{comments}
	Given an associative algebra $\mathcal A$, the natural question is to determine all commutators of $\mathcal A$. Shoda \cite{Shoda} proved that a matrix $C \in \mathbb M_n(F)$ is a commutator if and only if the trace of $C$ is zero. Wintner \cite{Wintner} proved that the identity in a unital Banach algebra is not a commutator. By passing to the Calkin algebra, Wintner's result immediately implies that a bounded operator on a Banach space which is of the form $\lambda I+K$ for some nonzero scalar $\lambda$ and a compact operator $K$ is not a commutator. Henceforth, researchers tried to characterize which operators on a given Banach space are commutators. The complete characterization of commutators in the Banach algebra $\mathcal B(\mathcal H)$ of all bounded operators on an infinite-dimensional Hilbert space $\mathcal H$ is due to Brown and Pearcy \cite{Brown:65}. They proved that a bounded operator $C$ on $\mathcal H$ is a commutator if and only if it is not of the form $\lambda I+K$ for some nonzero scalar $\lambda$ and some operator $K$ from the unique maximal ideal in $\mathcal B(\mathcal H)$. Apostol (\cite{Apostol:72, Apostol:73}) proved that a bounded operator on either $\ell^p$ ($1<p<\infty$) or $c_0$ is a commutator if and only if it is not of the form $\lambda I+K$ where $\lambda\neq 0$ and $K$ is compact. In the case of the Banach space $\ell^1$ the same characterization was obtained by Dosev in \cite{Dosev:09}. In the case of the Banach space, $\ell^\infty$ Dosev and Johnson \cite{DJ:10} proved that a bounded operator is a commutator if and only if it is not of the form $\lambda I+K$ where $\lambda\neq 0$ and $K$ is strictly singular.
	
	The study of positive commutators of positive operators on a given Banach lattice was initiated in \cite{Bracic:10}.
	The assumption on the positivity of $A, B$ and $C:= AB-BA$ may lead to some restrictions on the commutator. Namely, the authors proved that the positive commutator of positive compact operators is quasi nilpotent. They also posed a question of whether the same is true under the assumption that one of the operators is compact. This question was affirmatively and independently solved by R. Drnov\v sek \cite{Drn12} and N. Gao \cite{Gao14}. Inspired by a result of Schneeberger \cite{Schneeberger:71} asserting that a compact operator acting on a separable $L^p$ space $(1\leq p<\infty)$ is a commutator, in \cite{DK19} authors prove that a positive compact operator acting on a separable $L^p$ space $(1\leq p<\infty)$ is a commutator of positive operators. In \cite{DK19} authors provide a technical condition under which the answer to the proposed problem is affirmative for positive operators on $\ell^p$ space $(1\leq p\leq \infty)$ satisfying this technical condition.
\end{comments}

\communicator{Marko Kandi\'c}


\begin{problem}[Characterizing Koopman groups on Hilbert spaces]
    Use the Perron Frobenius spectral theory of positive $C_0$-groups to characterize unitary $C_0$-groups on a Hilbert space $H$ that are unitarily isomorphic to a Koopman group on an $L^2$-space.
\end{problem}

\begin{source}
    Manuscript of Rainer Nagel provided during the OPSO 2021 conference.
\end{source}


\section{Eventual positivity of semigroups}
	
\begin{problem}[Stability of eventually positive semigroups]
	Let $M$ be a von Neumann algebra equipped with a normal, semi-finite, and faithful trace $\tau$ and $p\in [1,\infty)$. 
	Let $(e^{tA})_{t\in [0,\infty)}$ be an individually eventually positive $C_0$-semigroup on $L^p(M,\tau)$.
	
	Is $s(A)=\omega_0(A)$?
\end{problem}

\begin{comments}
	A $C_0$-semigroup $(e^{tA})_{t\in [0,\infty)}$ on an ordered Banach space $E$ is said to be \emph{individually} eventually positive if for each $0\leq f\in E$, there exists $t_0\geq 0$ such that $e^{tA}f\geq 0$ for all $t\geq t_0$. If the time $t_0$ can be chosen independently of the initial value $f$, then we call the semigroup \emph{uniformly eventually positive}.
	
	For the case $p=1,2$, the answer is positive as shown in \cite[Theorem~6.2.6]{Glueck2016} (see also \cite[Theorem~7.8]{DanersGlueckKennedy2016a}); note that the proofs given in the aforementioned references are for the commutative $L^p$-spaces but they can be adapted to the non-commutative setting.
	If the semigroup is uniformly eventually positive and $E$ is the classical $L^p$-space, then a positive answer is given by Vogt \cite{Vogt2021}. However, for the non-commutative $L^p$-spaces with $p\neq 1,2$, the answer is not known even for positive semigroups.
	
	This problem was posed by Ralph Chill.
\end{comments}

\communicator{Sahiba Arora} 


\begin{problem}[Characterisation of eventual invariance by means of form methods]
    Let $(e^{tA})_{t \in [0,\infty)}$ be a $C_0$-semigroup on a Hilbert space $H$, where the operator $-A$ is associated to a form $a: V \times V \to H$ with form domain $V \subseteq H$. 
    Let $C \subseteq H$ be closed and convex.
    Give a characterisation, in terms of $a$ and $V$, of one or both of the following properties:
    \begin{enumerate}[(a)]
        \item 
        The semigroup leaves $C$ \emph{individually eventually invariant}, i.e., for each $c \in C$ exists a time $t_0 \in [0,\infty)$ such that $e^{tA}c \in C$ for all $t \ge t_0$.
    
        \item 
        The semigroup leaves $C$ \emph{uniformly eventually invariant}, i.e., there exists a time $t_0 \in [0,\infty)$ such that $e^{tA}c \in C$ for all for each $c \in C$ and all $t \ge t_0$.
    \end{enumerate}
    An answer would be particularly interesting in the case where $H$ is an $L^2$-space and $C$ is the usual positive cone in $H$ -- in this case, the properties mentioned in~(a) and~(b) become the properties \emph{individual eventual positivity} and \emph{uniform eventual positivity} that are mentioned in several further problems here.
\end{problem}

\begin{comments}
    \emph{Invariance} rather than \emph{eventual invariance} can indeed be characterised by means of $a$ and $V$.
    This is a very useful result due to Ouhabaz \cite{Ouhabaz1996}. The criterion becomes particularly simple if $H$ is an $L^2$-space and $C$ is the usual positive cone. If, in addition, $a$ is symmetric, this result is much older and goes back to Beurling and Deny \cite{bd59}.
    
  A characterization of (individual or uniform) eventual positivity by means of form methods might have the potential to considerably extent the applicability of the theory of eventual positivity.
\end{comments}

\communicator{Jochen Glück}


\section{$L^\infty$-bounds}

\begin{problem}[$L^\infty$-boundedness problem]\label{pb:linfty-bddness}
	Let $(\Omega,\mu)$ be a $\sigma$-finite measure space. Let $V$ be a dense subspace of $L^2(\Omega)$ 
	(complex-valued functions) and suppose that $\mathfrak{a}:V\times V\to\mathbb{C}$ is a closed sectorial 
	sesquilinar form, i.e.\ linear in the first and antilinear in the second argument. Here, sectoriality means that,
	for some $\theta\in[0,\frac{\pi}2)$,
$$
 \forall u\in V: \mathfrak{a}(u,u)\in\Sigma_\theta:=\{z\in\mathbb{C}\setminus\{0\}:|\mathrm{arg}\,z|\le\theta\}\cup\{0\},
$$
	and closedness of $\mathfrak{a}$ means that $V$ is complete for the norm 
$$
 \|u\|_{\mathfrak{a}}:=\left(\mathrm{Re}\,\mathfrak{a}(u,u)+\|u\|^2_{L^2(\Omega)}\right)^{1/2}.
$$
Let the linear operator $A$ in $L^2(\Omega)$ be associated with a $\mathfrak{a}$
in the sense that, for $u,h\in L^2(\Omega)$,
$$
 u\in D(A)\ \mbox{and}\ Au=h
 \quad\Longleftrightarrow\quad 
 u\in V\ \mbox{and}\ \forall v\in V: \mathfrak{a}(u,v)=\langle u,h\rangle,
$$
	where $\langle u,h\rangle:=\int_\Omega f\overline{h}\,d\mu$ denotes the usual scalar product in 
	$L^2(\Omega)$. Then $-A$ is negative generator of a bounded analytic semigroup $(T(\cdot))$ in $L^2(\Omega)$,
	which is contractive on the sector $\Sigma_{\frac{\pi}2-\theta}$.
	
	The semigroup $(T(t))_{t\ge0}$ is called {\em $L^\infty$-bounded}, if there exists $M>0$ such that 
\begin{equation}\label{eq:linfty-bdd}
 \|T(t)f\|_\infty\le M\|f\|_\infty,\qquad\mbox{for all $f\in L^2(\Omega)\cap L^\infty(\Omega)$.}
\end{equation}
	Can $L^\infty$-boundedness of $(T(t))_{t\ge0}$ be characterized in terms of the sesquilinear 
	form $\mathfrak{a}$?
\end{problem}

\begin{comments}
	The problem came up in discussions with S\"onke Blunck (at the end of the 1990s).	
	The case $M=1$ in \eqref{eq:linfty-bdd}, i.e., \emph{$L^\infty$-contractivity} of $(T(t))_{t\ge0}$, is characterized 
	in terms of the sesquilinear form $\mathfrak{a}$ by the well-know Beurling-Deny criterion (see, e.g., 
	\cite[Theorem~2.7]{Ouhabaz1992} or \cite[Section 2.2]{Ouhabaz2005}), namely by the condition
\begin{equation}\label{eq:BD-linfty}
  \forall u\in V: \mathrm{sgn}\,u(|u|-1)_+\in V\ \mbox{and}\ \mathrm{Re}\,\mathfrak{a}(u,\mathrm{sgn}\,u(|u|-1)_+)\ge0.
\end{equation}
	Here $v_+:=-((-v)\wedge 0)$, where $\wedge$ denotes the pointwise minimum, and 
      $\mathrm{sgn}\,u:=\frac{u}{|u|}1_{\{u\neq0\}}$ denotes the sign of the function $u$.

	A characterization of \emph{$L^\infty$-boundedness} of 	$(T(t))_{t\ge0}$ could certainly be 
	very useful (for the restrictions that $L^\infty$-contractivity imposes on the coefficients of second-order elliptic operators
	on domains we refer to \cite[Section 4.3]{Ouhabaz2005}). 

	It is clear that $(T(t))_{t\ge0}$ is $L^\infty$-bounded if there exists a function $g\in L^\infty(\Omega)$ with $g>0$ 
	$\mu$-a.e. and $1/g\in L^\infty(\Omega)$ satisfying
\begin{equation}\label{eq:linfty-contr-g}
 \|gT(t)f\|_\infty\le\|gf\|_\infty\qquad\mbox{for all $f\in L^\infty(\Omega)\cap L^2(\Omega)$,}
\end{equation}
	since the norm $f\mapsto\|gf\|_\infty$ is equivalent to $\|\cdot\|_\infty$. A modification of the condition 
	\eqref{eq:BD-linfty} characterizes \eqref{eq:linfty-contr-g}, namely
\begin{equation}\label{eq:BD-linfty-g}
 \forall u\in V: \mathrm{sgn}\,u(|u|-g)_+\in V\ \mbox{and}\ \mathrm{Re}\,\mathfrak{a}(u,\mathrm{sgn}\,u(|u|-g)_+)\ge0.
\end{equation}
	The proof can be done similar to the proof of equivalence of \eqref{eq:BD-linfty} and $L^\infty$-contractivity. 
	One can also resort to invariance results in $L^2(\Omega)$ and consider the closed convex set
$$
 K_g:=\{f\in L^2(\Omega):|f|\le g\ \mbox{$\mu$-a.e.}\}.
$$
	Then it is not hard to check that the projection $P_g$ of $L^2(\Omega)$ onto $K$  is given by 
	$P_gf:=\mathrm{sgn}\,f (|f|\wedge g)$ ($P_g(f)$ is the best approximation of $f$ in $K$).
	Then observe $u-P_g(u)=\mathrm{sgn}\,u(|u|-g)_+$ and apply \cite[Theorem 2.1]{Ouhabaz1996} (we also refer to
	\cite[Section 2.1]{Ouhabaz2005}). 
\end{comments}


In view of the comments to Problem~\ref{pb:linfty-bddness}, the following seems natural to ask.

\begin{problem}[$L^\infty$-contraction for a weight problem]\label{pb:linfty-contr-g}
	In the setting of Problem~\ref{pb:linfty-bddness} assume that the semigroup $(T(t))_{t\ge0}$ is $L^\infty$-bounded,
	i.e., satisfies \eqref{eq:linfty-bdd} for some $M\ge1$. Does there exist a function $g\in L^\infty(\Omega)$ satisfying 
	$g>0$ $\mu$-a.e., $1/g\in L^\infty(\Omega)$, and \eqref{eq:linfty-contr-g}?
\end{problem}

\begin{comments}
	 In general, this might be more than one can hope for.
\end{comments}

Hence we are led to the following.

\begin{problem}[Characterization of $L^\infty$-contraction for a weight]\label{pb:linfty-contr-g-char}
      In the setting of Problem~\ref{pb:linfty-bddness}, can one characterize those $L^\infty$-bounded semigroups 
	$(T(t))_{t\ge0}$, for which a function $g\in L^\infty(\Omega)$ satisfying $g>0$ $\mu$-a.e., $1/g\in L^\infty(\Omega)$,
	and \eqref{eq:linfty-contr-g} exists?
\end{problem}

\begin{comments}
	The question is whether $(T(t))_{t\ge0}$ can be made to be a contractive semigroup on the $\|\cdot\|_\infty$-closure of
	$L^2(\Omega)\cap L^\infty(\Omega)$ in $L^\infty(\Omega)$ for an equivalent norm of the special form
	$f\mapsto\|gf\|_\infty$. 
	There is, of course, and well-known from semigroup theory, a norm $\||\cdot\||$ on 
	$L^2(\Omega)\cap L^\infty(\Omega)$, equivalent to $\|\cdot\|_\infty$, such that
$$
 \||T(t)f\||\le\||f\||\qquad\mbox{for all $f\in L^\infty(\Omega)\cap L^2(\Omega)$.}
$$
	One can take $\||f\||:=\sup_{t\ge0}\|T(t)f\|_\infty$, which satisfies 
$$
 \|f\|_\infty\le\||f\||\le M\|f\|_\infty
$$
	by assumption \eqref{eq:linfty-bdd}. 
	So the problem might be seen as an $L^\infty$-counterpart to the question if every bounded 
	$C_0$-semigroup on a Hilbert space can be made contractive for an equivalent scalar product. The answer to this 
	question is known to be negative and there is a nice characterization of those bounded analytic $C_0$-semigroups that are 
	contractive for an equivalent scalar product in terms of bounded imaginary powers (and bounded $H^\infty$-calculus) 
	for the negative generator (see 	\cite{lemerdy}).  
\end{comments}


Still another question seems natural. 

\begin{problem}[Weight construction for $L^\infty$-contraction]\label{pb:linfty-contr-g-construct} 
	In the setting of Problem~\ref{pb:linfty-bddness}, assume that a function $g\in L^\infty(\Omega)$ satisfying $g>0$ 
	$\mu$-a.e., $1/g\in L^\infty(\Omega)$, and \eqref{eq:linfty-contr-g} exists. How can we find or construct such a 
	function $g$?
\end{problem}

\begin{comments}
	In the general situation of Problems~\ref{pb:linfty-bddness}, \ref{pb:linfty-contr-g}, \ref{pb:linfty-contr-g-char}, and 
	\ref{pb:linfty-contr-g-construct} it might well be that positivity 
	of the semigroup $(T(t))_{t\ge0}$ can help, i.e.\ the assumption that
	$f\ge0$ a.e. on $\Omega$ implies $T(t)f\ge0$ a.e. on $\Omega$ for all $t>0$. Recall that positivity of the semigroup 
	can be characterized in terms of the sesquilinear form $\mathfrak{a}$ (see \cite{Ouhabaz1992}, \cite{Ouhabaz2005}).

	Assume for the following that the semigroup $(T(t))_{t\ge0}$ is positive and that the measure space 
	$(\Omega,\mu)$ is \emph{finite}. Then $L^\infty(\Omega)\subseteq L^2(\Omega)$
	and hence $T(t)f\in L^2(\Omega)$ is defined for any $f\in L^\infty(\Omega)$. In this situation, $L^\infty$-contractivity 
	of $(T(t))_{t\ge0}$ is characterized by $T(t)1_\Omega\le 1_\Omega$ $\mu$-a.e. for all $t>0$ 
	where $1_\Omega$ denotes the characteristic function of $\Omega$. 

	Now let $g\in L^\infty(\Omega)$ such that 
	$g>0$ $\mu$-a.e. on and 	$1/g\in L^\infty(\Omega)$. 
	Then, consequently, \eqref{eq:linfty-contr-g} is characterized by $T(t)g^{-1}\le g^{-1}$ $\mu$-a.e. for all $t>0$ 
	(since this is equivalent to $L^\infty$-contractivity of the positive semigroup $(gT(t)g^{-1})_{t\ge0}$). 
	
	In particular, one has \eqref{eq:linfty-contr-g} for $g\in L^\infty(\Omega)$ with $g>0$ $\mu$-a.e. if 
	$h:=1/g\in L^\infty(\Omega)$ is an \emph{eigenfunction} for an eigenvalue $\lambda\ge0$ of $A$: Recalling that
	$-A$ is the generator of $(T(t))_{t\ge0}$ we obtain $T(t)h= e^{-\lambda t}h\le h$ $\mu$-a.e. for all $t>0$.

	Specializing further, take $\Omega\subset\mathbb{R}^d$ a bounded domain (with $\mu$ the Lebesgue measure) and 
	the usual Dirichlet form $\mathfrak{a}(u,v):= \int_\Omega \nabla u\cdot\overline{\nabla v}\,dx$ with form domain 
	$V:=V_N:=H^1(\Omega)$ (Neumann boundary conditions) or $V:=V_D:=H^1_0(\Omega)$ (Dirichlet boundary 	
	conditions), and denote the associated operators by $A_N$ and $A_D$, respectively (the negative Laplacian on 
	$\Omega$ with Neumann/Dirichlet boundary conditions). 

	Both semigroups are well-known 
	to be positive and $L^\infty$-contractive, i.e.\ they satisfy \eqref{eq:linfty-contr-g} for $g=1_\Omega$. Now $1_\Omega$
	is an eigenfunction of $A_N$ for the eigenvalue $0$. However, for $A_D$ we do not even have $1_\Omega\in V_D$. We 
	do have an eigenfunction $h\in L^\infty(\Omega)$, $h>0$ $\mu$-a.e. on $\Omega$, for the first eigenvalue 
	$\lambda_0>0$ of $A_D$, but $1/h\not\in L^\infty(\Omega)$ due to $h\in V_D=H^1_0(\Omega)$. Hence, considering 	
	(positive) eigenfunctions are not sufficient, in general.

	It might be more adequate to consider \emph{positive subeigenfunctions}, we refer to \cite[II-C Section  3]{Arendt1986}
	for the notion of positive sub eigenvectors and their role in the characterization of the positivity of semigroups on Banach lattices.
\end{comments}

\communicator{Peer Kunstmann}

\section{Further problems}


\begin{problem}[Closedness of ranges]
    Consider the following properties of the generation $A$ of a $C_0$-semigroups $(T(t))_{t\geq0}$.
    \begin{enumerate}[(i)]
        \item
        $\operatorname{rg} A$ closed.
        
        \item 
        $\operatorname{rg}(1-T(t))$ closed for one $t>0$.
    \end{enumerate}
    Does assertion (ii) imply assertion (i)?
\end{problem}

\textit{Comment.} This is a version of a "spectral inclusion theorem". The converse implication does not hold.

\medskip

\begin{source}
	\bari{Yuri Latushkin}
\end{source}


\begin{problem}[Lower bounds]
    Let $(T(t))_{t \geq 0}$ be a contraction semigroup with generator $A$ on a Hilbert space $H$. Consider the following properties:
    \begin{enumerate}[(i)]
        \item
        There exists $m \geq 0$ such that
        \begin{align*}
            \|(\lambda - A)x\| \geq m |\operatorname{Re} \lambda| \|x\|
            \quad \text{for all } \operatorname{Re} \lambda < 0 
            \text{ and all } x \in H.
        \end{align*}
        
        \item
        There exists $m_1 > 0$ such that
        \begin{align*}
            \|T(t)x\| \geq m_1 \|x\|
            \quad \text{for all } t \geq 0, \, x \in H.
        \end{align*}
    \end{enumerate}
    Does assertion $(i)$ imply assertion $(ii)$?
\end{problem}

\textit{Comment.} Condition $(ii)$ always implies $(i)$. This is not true if $(T(t))_{t \geq 0}$ is only bounded, but it holds if $\lambda \in \rho(A)$ with $\operatorname{Re} \lambda < 0$.
\medskip 

\textit{Comment} (Charles Batty, April 2021): 
There is no answer yet for this open question. In a paper by Xu and Shang \cite{Xu2009}, a related result on Banach spaces is stated in Theorem 2.4. Their proof is seriously flawed, but a correct proof is given in a paper by Geyer and Batty \cite[Theorem 5.4, Example 5.6]{BaGe2017}.

\medskip

\begin{source}
	\bari{Birgit Jacob and Hans Zwart}
\end{source}


\begin{problem}[Backward uniqueness property]
    Study the "backward uniqueness property", i.e., characterize injective $C_0$-semigroups. Apply the results to the backward uniqueness property for non-autonomous Cauchy problems $u'(t) = A(t)u(t)$, where $A(t)$ is sectorial, by examining the corresponding evolution semigroup. See Lunardi's work on evolution equations in Banach spaces \cite{lunardi2009}.
\end{problem}


\begin{source}
    \bari{Alessandra Lunardi}
\end{source}

 



\begin{problem}[Rational calculus and contractivity]
	~
	\begin{enumerate}[(i)]
   		\item Does every bounded $C_0$-semigroup on a Hilbert space have a bounded rational calculus?
    		\item When is a $C_0$-semigroups on a Hilbert space similar to a contraction semigroup?
	\end{enumerate}
\end{problem}

\begin{comments}
   The answer to the first question is no. A bounded semigroup on a Hilbert space has a bounded rational calculus if and only if it has a bounded $H^{\infty}$-calculus. For more details, see \cite{Haase2003}. 
\end{comments}

\begin{comments}(Charles Batty)
	In Problem (i) and the Comment above, the "boundedness" of the rational calculus is being interpreted as meaning boundedness with respect to the $H^\infty$-norm.   There are one or two issues as to what is the domain of those functions and whether or not the generator is injective, but the answer is correct. See also  \cite[Sections 5.3.4, section  5.3.5] {Haase2006} for more details.

	Instead of considering the $H^\infty$-norm, one may consider Banach algebras in different norms that are embedded in $H^\infty(\mathbb{C}_+)$, where $\mathbb{C}_+$ is the open right half-plane.   One example is the (Hille)-Phillips calculus, where the norm comes from measures on $[0,\infty)$.   Alexander Gomilko, Yuri Tomilov and Hans Zwart \cite{GomilkoZwart2007} have shown that if $-A$ is the generator of a bounded $C_0$-semigroup on a Hilbert space, then there is a bounded $\mathcal{B}$-calculus for $A$ where   $\mathcal{B}$ is a Banach algebra of "analytic Besov" functions on the right half-plane. This algebra is considerably bigger than the Phillips algebra, and the $\mathcal{B}$-norm is considerably smaller than the Phillips norm but it is bigger than the $H^\infty$-norm.

	On Banach spaces, the $\mathcal{B}$-calculus exists if and only if $A$ satisfies the condition introduced by Gomilko, and independently by Shi and Feng, in 1999 and 2000.   In particular, it exists if  $A$ is sectorial of angle less than $\pi/2$, so $-A$ generates a bounded holomorphic $C_0$-semigroup.   For those operators, there are two further calculi, $\mathcal{D}$-calculus and $\mathcal{H}$-calculus, which extend the calculus to larger classes of functions than $\mathcal{B}$, with smaller norms.  
\end{comments}

\begin{relevance}
	This problem is related to the inverse generator problem and questions concerning the powers of the co-generator of bounded semigroups. For both problems, the answers have been negative in general, and some positive partial answers have been obtained, but the answer for bounded semigroups on Hilbert spaces is unknown.  The extended calculi provide systematic ways to approach such problems, instead of using ad hoc methods each time.  For the functions $((z-1)(z+1)^{-1})^n$  the Hille-Phillips norm grows like $n^{1/2}$, the $\mathcal{B}$-norm grows like $\log n$, and the $\mathcal{D}$-norms are uniformly bounded. 
\end{relevance}

\begin{source}
    \bari{Hans Zwart}
\end{source}


\begin{problem}[A.E. Teretenkov]
    Let $ \mathcal{B} $ be a Banach space. Let $ \mathcal{P} $  be a projection on a \textit{finite-dimensional} Banach subspace of $ \mathcal{B} $. Let $ \mathcal{L}^0 + \lambda \mathcal{L} $ be a generator of  $ C_0 $-semigroup $ \mathcal{U}_t^{\lambda} $ on $ \mathcal{B} $ for all $  \lambda \in [0, \lambda_{\rm sup})$. Let $ \mathcal{U}_t^{0} $  leave both $ \mathcal{P} \mathcal{B} $ and its complement  $ (I- \mathcal{P} )\mathcal{B} $ invariant, let $ \mathcal{P} \mathcal{L} \mathcal{P} = 0 $. Denote $ \mathcal{L} _t \equiv (\mathcal{U}_{t}^{0})^{-1} \mathcal{L} \mathcal{U}_t^{0}  $. Let the integrals
\begin{equation*}
	 \int_{-\infty}^t dt_1 \ldots \int_{-\infty}^{t_{k-1}} dt_k \mathcal{P} \mathcal{L}_{t_1}  \ldots   \mathcal{L}_{t_k} \mathcal{P}, \qquad k = 1, \ldots, n+1
\end{equation*}
be finite for all $ t \geqslant 0 $. If is it possible to find such a $ \lambda $-dependent operator $ r^{n, \lambda}  $ on $ \mathcal{P}\mathcal{B} $ , which is polynomial in $ \lambda $, and $ \lambda $-dependent  semigroup $ u_t^{n, \lambda} $ on $ \mathcal{P}\mathcal{B} $, whose generator is polynomial in $ \lambda $,   such that 
\begin{equation*}
	\mathcal{P}( \mathcal{U}_{\frac{t}{\lambda^2}}^0)^{-1} \mathcal{U}_{\frac{t}{\lambda^2}}^{\lambda} \mathcal{P} =  u_t^{n, \lambda} r^{n, \lambda}   + O(\lambda^{2n+2}), \qquad \lambda \rightarrow + 0
\end{equation*}
for all $ t>0 $ ?

\medskip

Let us emphasize that we do not assume here asymptotic behaviour to be uniform in $ t $, it is just assumed for each fixed $ t > 0 $. If there are counterexamples, what further restrictions should be assumed to obtain such an asymptotic estimate? 
\end{problem}

\begin{relevance}
This problem is important for the derivation of perturbative corrections to Markovian quantum master equations. It seems to be a possible direction of generalization of the classical results by E.B.~Davies \cite{Davies1974}  to higher orders of perturbation theory in $ \lambda $ and seems to be held in a simple example discussed in \cite{Teretenkov2020}. It also seems to be necessary for strict perturbative derivation of master equations recently obtained by A.S.~Trushechkin \cite{Trushechkin2021}.
\end{relevance}

\begin{problem}[Super-fast converging Chernoff approximations example] Provide a meaningful example of Chernoff approximations that converge to a $C_0$-semigroup at least exponentially or prove that it is not possible in a meaningful setting.

\textit{Informal statement.} Provide a concrete example of a super-fast (at least exponentially) converging Chernoff approximation to a semigroup generated by a first-order differential operator with variable coefficients, with Chernoff function explicitly expressed in terms of these coefficients and without circulus vitiosus of any kind, e.g.\ in the form of using the values of the semigroup for creating a Chernoff function for it.

\textit{Formal statement.} Consider Banach space $UC_b(\mathbb{R},\mathbb{R})$ of all bounded and uniformly continuous real-valued functions of one real variable, with the uniform norm $\|f\|=\sup_{x\in\mathbb{R}}|f(x)|$. Recall that the space $C_b^\infty(\mathbb{R},\mathbb{R})$ of all functions that are bounded and have bounded derivatives of all orders is a dense linear subspace in $UC_b(\mathbb{R},\mathbb{R})$. Suppose that arbitrary functions $a,b,c\in UC_b(\mathbb{R},\mathbb{R})$ satisfying $a(x)\geq a_0\equiv\mathrm{const}>0$ for all $x\in\mathbb{R}$ are given. Define an operator $L$ by the equality 
\[
(Lf)(x)=a(x)f''(x)+b(x)f'(x)+c(x)f(x) \text{ for all }  f\in C_b^\infty(\mathbb{R},\mathbb{R}), x\in\mathbb{R}
\] and recall that the closure of $L$ generates a $C_0$-semigroup $(e^{tL})_{t\geq 0}$ in $UC_b(\mathbb{R},\mathbb{R})$, satisfying 
\[
\|e^{tL}\|\leq e^{t\max(0, \sup\limits_{x\in\mathbb{R}} c(x))} \quad \text{ for all } \quad  t\geq 0.
\]
 The problem is to find a family $(S(t))_{t\geq 0}$ of everywhere defined linear bounded operators in $UC_b(\mathbb{R},\mathbb{R})$ which has the following three properties:

\begin{itemize}
\item $S$ and $L$ satisfy the conditions of at least one version of the Chernoff theorem on approximation of operator semigroups, i.e.\ $S$ is Chernoff-equivalent to the semigroup $(e^{tL})_{t\geq 0}$, i.e.\ $S$ is a Chernoff function for the operator $L$;

\item There exists a dense linear subspace $D$ in $UC_b(\mathbb{R},\mathbb{R})$ and constants $M>0, q>1, n_0>0$ such that 
\[
\|e^{tL}f-S(t/n)^nf\|\leq M/q^{n} \quad  \text{ for all }  \quad n\geq n_0 \text{ and all } f\in D;
\]

\item $S$ is defined by an explicitly given formula $(S(t)f)(x)=F(a,b,c,f,x)$ and the calculation of $e^{tL}$ is not included into this formula, even in a hidden way, so expression $e^{tL}=\lim\limits_{n\to\infty}S(t/n)^n$ is free from circulus vitiosus of any kind and can be used as a practical method of calculation of $e^{tL}$.

\end{itemize}
\end{problem}

\begin{comments}
After the Chernoff theorem was published in 1968, it has found many applications in random processes theory, partial differential equations, functional integration and other fields closely related to operator semigroups theory, see \cite{Butko2020}. 
Chernoff's original theorem states the fact of convergence however says nothing about the rate of convergence. It was shown (see papers \cite{Zagrebnov2020, Zagrebnov2022, Gomilko2019, CACHIA2001176, IcTa1997, NeZa1998}) that, as $n$ tends to infinity, in some cases $\|e^{tL}-S(t/n)^n\|$ tends to zero with a rate of $const/n$, and under some conditions may reach $const/n^2$. However the rate of tending of $\|e^{tL}f-S(t/n)^nf\|$ to zero depends on $f$ heavily (see \cite{Prudnikov2020}, \cite{Dragunova2023}) and may be slower than $const/n$. Model examples show that $\|e^{tL}f-S(t/n)^nf\|$ can tend to zero with arbitrary high or arbitrary low rate and this is also possible in a general setting -- for arbitrary never-zero semigroup (see  \cite{Galkin2022}). Summing up, we see that the rate of convergence of Chernoff approximations depends on all three elements: $e^{tL}$, $S$, $f$. 

If $e^{tL}$ and $S$ has the same Taylor's polynomial of order $k$ and $S$ is close to its Taylor polynomial in some sense (see  \cite{Galkin2022}), and $f$ is good enough related to $L$, then $\|e^{tL}f-S(t/n)^nf\|\leq const/n^k$. If the rate of convergence of Chernoff approximations is higher than $const/n$ then we say that they are fast converging Chernoff approximations, and if the rate is higher than $const/n^k$ for all $k\in\mathbb{N}$ then we say that they are super fast converging Chernoff approximations. It is expected that Chernoff functions $S$ which provide super fast converging approximations should satisfy $\frac{d^k}{dt^k}S(t)\big|_{t=0}=\frac{d^k}{dt^k}e^{tL}\big|_{t=0}$ for each $k\in\mathbb{N}$. 

Examples of fast converging approximations for the operator $L$ given by $Lf=af''+bf'+cf$ were constructed for $k=2$ and it is clear (see  \cite{Vedenin2020})how to extend this method to arbitrary positive integer values of $k$
involves a significant amount of work to be done. If we already know the semigroup $e^{tL}$ then we can construct Chernoff function $S$ that provides Chernoff approximations $S(t/n)^n$ that converge to $e^{tL}$ with arbitrarily high rate (even with infinitely high rate if we take $S(t)=e^{tL}$) but they are not useful in practice because we employ the semigroup to approximate it falling into circulus vitiosus. Effective methods for constructing super fast converging Chernoff approximations for this operator $L$ are not known at present. However, there are no fundamental objections to their existence either, which presents a challenge to either find such approximations or to prove that they are impossible for some reason.
\end{comments}

\begin{relevance}
Examples of fast-converging Chernoff approximations known at present all involve exponential growth in computational complexity: to calculate $S(t/n)^n$ we need to perform $const\cdot p^n$ arithmetic operations where $p$ ranges from 3 to 9 for different variants of $S$ known at present. Even worse is the fact that large values of $k$ in the estimate $\|e^{tL}f-S(t/n)^nf\|\leq const/n^k$ result in large values of $p$. However if one can find Chernoff function $S$ that provides $\|e^{tL}f-S(t/n)^nf\|\leq const/q^{n}$ for $q>1$ then this exponential decay of the approximation error may balance the (hoped to be still not worse than exponential) growth of computational complexity. If it happens then these super fast converging Chernoff approximations can be a good alternative to existing numerical methods of solving various linear differential equations with variable coefficients: evolution PDEs (Schr\"odinger type equations with complex coefficients, parabolic equations with real coefficients) as well as ODEs and elliptic PDEs. 
\end{relevance}

\communicator{Ivan D. Remizov}

\begin{acknowledgement}
The authors are grateful to everyone who has communicated open problems. Their contributions have been invaluable in advancing our research. A large portion of problems was found in the unpublushed list of problems provided by Rainer Nagel, a large portion of editorial work was done by Jochen Gl\"uck, the authors are thankful to both of them. 

If you have new problems that deserve being included to the list or you have solutions or comments to existing problems please email editors of the list. The latest version of this list is here:  https://arxiv.org/abs/2410.00416.
\end{acknowledgement}

{}

\end{document}